\documentclass[10pt,twoside]{article}

\usepackage{amsmath,enumerate,amsfonts,amssymb,amsthm}
\evensidemargin\oddsidemargin \itemsep=0pt
\theoremstyle{plain}

\def\fnote#1{\footnote}

\input amssym.def
\input amssym

\begin{document}

\pagestyle{myheadings} \markboth {{\small\sc Constant Jacobi
osculating rank of $U(3)/(U(1) \times U(1) \times U(1))$ -
Appendix -}} {{\small\sc T. Arias-Marco}}

\thispagestyle{empty}

\begin{center}

{\Large\bf Constant Jacobi osculating rank of $\mathbf{U(3)/(U(1)
\times U(1) \times U(1))}$ \\\vspace{0.2cm}

- Appendix -} \vskip 5mm

%%%%%%%%%%%%%%%%%%%%%%%%%%%%%%%%%%%%%%%%%%%%%%%%%%%%%%%%%%%%%%%%
%
% Author(s), affiliation(s) and email(s). For instance:
%
%%%%%%%%%%%%%%%%%%%%%%%%%%%%%%%%%%%%%%%%%%%%%%%%%%%%%%%%%%%%%%%%

{\large\bf Teresa Arias-Marco}\footnote{The author's work has been
supported by D.G.I.~(Spain) and FEDER Project MTM 2007-65852, the
network MTM2008-01013-E/ and by DFG Sonderforschungsbereich~647.}

\vskip 5mm {\it Departamento de Matem\'aticas,
Universidad de Extremadura, \\
Avenida de Elvas s/n, 06071 Badajoz, Spain.\\
email: ariasmarco@unex.es}
\end{center}
\begin{center}
{\small\it Dedicated to Salud Bartoll.}
\end{center}
%\medskip
\begin{abstract}
This is the appendix of the paper [T. Arias-Marco, Constant Jacobi
osculating rank of $U(3)/(U(1) \times U(1) \times U(1))$,  Arch.
Math. (Brno) 45 (2009), 241--254] where we obtain an interesting
relation between the covariant derivatives of the Jacobi operator
valid for all geodesic on the flag manifold $M^6=U(3)/(U(1) \times
U(1) \times U(1))$. As a consequence, an explicit expression of
the Jacobi operator independent of the geodesic can be obtained on
such a manifold. Moreover, we show the way to calculate the Jacobi
vector fields on this manifold by a new formula valid on every
g.o. space.

\vspace{1\baselineskip}


\end{abstract}

\vspace{1\baselineskip}

\par
\medskip

\section*{Appendix}

$\mathcal{J}_0^{2)} = \left(\mathcal{J}_{ij}^2(0)\right)$, $i, j =
1, \ldots, 6$, where
\begin{equation}\label{J2}
\begin{array}{l}
\mathcal{J}_{11}^2(0) =  - \frac{3}{4} \left( x_4^2 (x_2^2 + x_3^2
+ x_5^2 + x_6^2)
- 2 (x_2 x_6 + x_3 x_5)^2 \right),  \\

\mathcal{J}_{12}^2(0) =  \frac{3}{8} \left(x_4 x_5 (x_1^2 + x_2^2
+ x_4^2 + x_5^2) -
2 (x_2 x_6 + x_3 x_5) (x_1 x_6 + x_3 x_4) \right), \\

\mathcal{J}_{13}^2(0) = \frac{3}{8} \left(x_4 x_6 (x_1^2 + x_3^2 +
x_4^2 + x_6^2) -
2 (x_2 x_6 + x_3 x_5) (x_1 x_5 + x_2 x_4) \right), \\

\mathcal{J}_{14}^2(0) =  \frac{3}{4} \left( x_1 x_4 (x_2^2 + x_3^2
+ x_5^2 +
x_6^2) + 2 (x_2 x_6 + x_3 x_5) (x_2 x_3 - x_5 x_6)  \right),  \\

\mathcal{J}_{15}^2(0) = - \frac{3}{8} \left(x_2 x_4 (x_1^2 + x_2^2
+ x_4^2 + x_5^2)
+ 2 (x_2 x_6 + x_3 x_5) (x_1 x_3 - x_4 x_6)) \right), \\

\mathcal{J}_{16}^2(0) = - \frac{3}{8} \left(x_3 x_4 (x_1^2 + x_3^2
+ x_4^2 + x_6^2)
+ 2 (x_2 x_6 + x_3 x_5) (x_1 x_2 - x_4 x_5) \right),\\

\\

\mathcal{J}_{22}^2(0) =  - \frac{3}{4} \left( x_5^2 (x_1^2 + x_3^2
+ x_4^2 + x_6^2)
- 2 (x_1 x_6 + x_3 x_4)^2 \right),  \\

\mathcal{J}_{23}^2(0) = \frac{3}{8} \left(x_5 x_6 (x_2^2 + x_3^2 +
x_5^2 + x_6^2) -
2 (x_1 x_6 + x_3 x_4) (x_1 x_5 + x_2 x_4) \right), \\

\mathcal{J}_{24}^2(0) = - \frac{3}{8} \left(x_1 x_5 (x_1^2 + x_2^2
+ x_4^2 + x_5^2)
+ 2 (x_1 x_6 + x_3 x_4) (x_2 x_3 - x_5 x_6)) \right), \\

\mathcal{J}_{25}^2(0) =  \frac{3}{4} \left( x_2 x_5 (x_1^2 + x_3^2
+ x_4^2 +
x_6^2) + 2 (x_1 x_6 + x_3 x_4) (x_1 x_3 - x_4 x_6)  \right),  \\

\mathcal{J}_{26}^2(0) = - \frac{3}{8} \left(x_3 x_5 (x_2^2 + x_3^2
+ x_5^2 + x_6^2)
+ 2 (x_1 x_6 + x_3 x_4) (x_1 x_2 - x_4 x_5) \right),\\

\\

\mathcal{J}_{33}^2(0) =  - \frac{3}{4} \left( x_6^2 (x_1^2 + x_2^2
+ x_4^2 + x_5^2)
- 2 (x_1 x_5 + x_2 x_4)^2 \right),  \\

\mathcal{J}_{34}^2(0) = - \frac{3}{8} \left(x_1 x_6 (x_1^2 + x_3^2
+ x_4^2 + x_6^2)
+ 2 (x_1 x_5 + x_2 x_4) (x_2 x_3 - x_5 x_6)) \right), \\

\mathcal{J}_{35}^2(0) = - \frac{3}{8} \left(x_2 x_6 (x_2^2 + x_3^2
+ x_5^2 + x_6^2)
+ 2 (x_1 x_5 + x_2 x_4) (x_1 x_3 - x_4 x_6) \right),\\

\mathcal{J}_{36}^2(0) =  \frac{3}{4} \left( x_3 x_6 (x_1^2 + x_2^2
+ x_4^2 +
x_5^2) + 2 (x_1 x_5 + x_2 x_4) (x_1 x_2 - x_4 x_5)  \right),  \\

\\

\mathcal{J}_{44}^2(0) =  - \frac{3}{4} \left( x_1^2 (x_2^2 + x_3^2
+ x_5^2 + x_6^2)
- 2 (x_2 x_3 - x_5 x_6)^2 \right),  \\

\mathcal{J}_{45}^2(0) = \frac{3}{8} \left(x_1 x_2 (x_1^2 + x_2^2 +
x_4^2 + x_5^2) -
2 (x_2 x_3 - x_5 x_6) (x_1 x_3 - x_4 x_6) \right), \\

\mathcal{J}_{46}^2(0) = \frac{3}{8} \left(x_1 x_3 (x_1^2 + x_3^2 +
x_4^2 + x_6^2) -
2 (x_2 x_3 - x_5 x_6) (x_1 x_2 - x_4 x_5) \right), \\

\\

\mathcal{J}_{55}^2(0) =  - \frac{3}{4} \left( x_2^2 (x_1^2 + x_3^2
+ x_4^2 + x_6^2)
- 2 (x_1 x_3 - x_4 x_6)^2 \right), \\

\mathcal{J}_{56}^2(0) = \frac{3}{8} \left(x_2 x_3 (x_2^2 + x_3^2 +
x_5^2 + x_6^2) -
2 (x_1 x_3 - x_4 x_6) (x_1 x_2 - x_4 x_5) \right), \\

\\

\mathcal{J}_{66}^2(0) =  - \frac{3}{4} \left( x_3^2 (x_1^2 + x_2^2
+ x_4^2 + x_5^2)
- 2 (x_1 x_2 - x_4 x_5)^2 \right).  \\

\end{array}
\end{equation}

\newpage

$\mathcal{J}_0^{3)} = \left(\mathcal{J}_{ij}^3(0)\right)$, $i, j =
1, \ldots, 6$, where
\begin{equation}\label{J3}
\begin{array}{ccl}

\mathcal{J}_{11}^3(0) & = &   \frac{9}{8\sqrt{2}} \ x_4 (x_2^2 -
x_3^2 + x_5^2 - x_6^2)(x_2 x_6 + x_3 x_5),   \\

\mathcal{J}_{12}^3(0) & = & \frac{3}{16\sqrt{2}} \left( x_3 (
x_1^2 (x_4^2 - 4 x_5^2) + x_2^2 (4 x_4^2 - x_5^2) + (x_4^2 -
x_5^2) (7 x_3^2 +
x_4^2 + x_5^2 + 7 x_6^2))  \right. \\
 &  & + \left. x_1 x_4 x_6 (x_1^2 + 4 x_2^2 + 7 x_3^2 + x_4^2 + 4 x_5^2 + 7 x_6^2) - x_2 x_5 x_6
 (4 x_1^2 + x_2^2 + 7 x_3^2 + 4 x_4^2 + x_5^2 + 7 x_6^2 )\right), \\

\mathcal{J}_{13}^3(0) & = & \frac{3}{16\sqrt{2}} \left( - x_4 (
x_1 x_5 + x_2 x_4 )
(x_1^2 + 7 x_2^2 + 4 x_3^2 + x_4^2 + 7 x_5^2 + 4 x_6^2))  \right. \\
 &  & + \left. x_6( x_2 x_6 +  x_3 x_5 ) (4 x_1^2 + 7 x_2^2 + x_3^2 + 4 x_4^2 +
7 x_5^2 + x_6^2 ) \right),\\

\mathcal{J}_{14}^3(0) & = &   \frac{9}{16\sqrt{2}} \ (x_2^2 -
x_3^2 + x_5^2 - x_6^2)
(x_1 (x_2 x_6 - x_3 x_5) + x_4 (x_2 x_3 - x_5 x_6)),  \\

\mathcal{J}_{15}^3(0) & = & \frac{3}{16\sqrt{2}} \left(  x_4 ( x_1
x_3 - x_4 x_6 )
(x_1^2 + 4 x_2^2 + 7 x_3^2 + x_4^2 + 4 x_5^2 + 7 x_6^2))  \right. \\
 &  & + \left. x_2( x_2 x_6 +  x_3 x_5 ) (4 x_1^2 + x_2^2 + 7 x_3^2 + 4 x_4^2 + x_5^2 + 7 x_6^2 ) \right),\\

 \mathcal{J}_{16}^3(0) & = & \frac{3}{16\sqrt{2}} \left( x_5 ( x_1^2 (x_4^2 - 4
x_3^2) + x_6^2 (4 x_4^2 - x_3^2) + (x_4^2 - x_3^2) (7 x_2^2 +
x_3^2 + x_4^2 + 7 x_5^2))  \right. \\
 &  & - \left.  x_1 x_2 x_4 (x_1^2 + 7 x_2^2 + 4 x_3^2 + x_4^2 + 7 x_5^2 + 4 x_6^2) - x_6 x_2
 x_3 (4 x_1^2 + 7 x_2^2 + x_3^2 + 4 x_4^2 + 7 x_5^2 + x_6^2 )\right), \\

 \\

 \mathcal{J}_{22}^3(0) & = &  \frac{9}{8\sqrt{2}} \ x_5 ( x_1^2 -
x_3^2 + x_4^2 - x_6^2)(x_1 x_6 + x_3 x_4),   \\

\mathcal{J}_{23}^3(0) & = & \frac{3}{16\sqrt{2}} \left(  x_5 ( x_1
x_5 + x_2 x_4 )
(7 x_1^2 + x_2^2 + 4 x_3^2 + 7 x_4^2 + x_5^2 + 4 x_6^2))  \right. \\
 &  & - \left. x_6( x_1 x_6 +  x_3 x_4 ) (7 x_1^2 + 4 x_2^2 + x_3^2 + 7 x_4^2 +
4 x_5^2 + x_6^2 ) \right),\\

\mathcal{J}_{24}^3(0) & = & \frac{3}{16\sqrt{2}} \left( x_6 (
x_2^2 ( 4 x_1^2  - x_5^2 ) + x_4^2 ( x_1^2 - 4 x_5^2) + (x_1^2 -
x_5^2) (x_1^2 + 7 x_3^2 +
 x_5^2 + 7 x_6^2))  \right. \\
 &  & - \left. x_2 x_3 x_5 (4 x_1^2 + x_2^2 + 7 x_3^2 + 4 x_4^2 + x_5^2 + 7 x_6^2 )- x_4 x_1 x_3
  (x_1^2 + 4 x_2^2 + 7 x_3^2 + x_4^2 + 4 x_5^2 + 7 x_6^2)  \right), \\

\mathcal{J}_{25}^3(0) & = &   \frac{9}{16\sqrt{2}} \ (x_1^2 -
x_3^2 + x_4^2 - x_6^2)
(x_2 (x_1 x_6 + x_3 x_4) - x_5 (x_1 x_3 + x_4 x_6)),  \\

\mathcal{J}_{26}^3(0) & = & \frac{3}{16\sqrt{2}} \left( x_4 (
x_2^2 (4 x_3^2 - x_5^2) + x_6^2 (x_3^2 - 4 x_5^2) + (x_3^2 -
x_5^2) (7 x_1^2 + x_3^2 +
7 x_4^2 + x_5^2))  \right. \\
 &  & + \left. x_2 x_1 x_5 (7 x_1^2 + x_2^2 + 4 x_3^2 + 7 x_4^2 +  x_5^2 + 4 x_6^2) + x_6
 x_1 x_3 (7 x_1^2 + 4 x_2^2 + x_3^2 + 7 x_4^2 + 4 x_5^2 + x_6^2 )\right), \\

\\

\mathcal{J}_{33}^3(0) & = &  \frac{9}{8\sqrt{2}} \ x_6 ( x_1^2 -
x_2^2 + x_4^2 - x_5^2)(x_1 x_5 + x_2 x_4),   \\

\mathcal{J}_{34}^3(0) & = & \frac{3}{16\sqrt{2}} \left(  x_6 ( x_2
x_3 - x_5 x_6 )
(4 x_1^2 + 7 x_2^2 + x_3^2 + 4 x_4^2 + 7 x_5^2 + x_6^2))  \right. \\
 &  & + \left. x_1( x_1 x_5 +  x_2 x_4 ) (x_1^2 + 7 x_2^2 + 4 x_3^2 + x_4^2 +
7 x_5^2 + 4 x_6^2 ) \right),\\

\mathcal{J}_{35}^3(0) & = & \frac{3}{16\sqrt{2}} \left(  x_6 (-
x_1 x_3 + x_4 x_6 )
(7 x_1^2 + 4 x_2^2 + x_3^2 + 7 x_4^2 + 4 x_5^2 + x_6^2))  \right. \\
 &  & - \left. x_2( x_1 x_5 +  x_2 x_4 ) (7 x_1^2 + x_2^2 + 4 x_3^2 + 7 x_4^2 +
x_5^2 + 4 x_6^2 ) \right),\\

\mathcal{J}_{36}^3(0) & = &   \frac{9}{16\sqrt{2}} \ (x_1^2 -
x_2^2 + x_4^2 - x_5^2)
(- x_3 (x_1 x_5 + x_2 x_4) + x_6 (x_1 x_2 - x_4 x_5)),  \\

 \\

 \mathcal{J}_{44}^3(0) & = &  \frac{9}{8\sqrt{2}} \ x_1 ( x_2^2 -
x_3^2 + x_5^2 - x_6^2)(- x_2 x_3 + x_5 x_6),   \\

\mathcal{J}_{45}^3(0) & = & \frac{3}{16\sqrt{2}} \left( x_3 (-
x_1^2 (x_4^2 + 4 x_5^2) + x_2^2 (4 x_4^2 + x_5^2) - (x_1^2 -
x_2^2) ( x_1^2 + x_2^2 +
7 x_3^2 + 7 x_6^2))  \right. \\
 &  & + \left. x_4 x_1 x_6 (x_1^2 + 4 x_2^2 + 7 x_3^2 + x_4^2 + 4 x_5^2 + 7 x_6^2) +
 x_5 x_2 x_6 (4 x_1^2 + x_2^2 + 7 x_3^2 + 4 x_4^2 +  x_5^2 + 7 x_6^2 )\right), \\

\mathcal{J}_{46}^3(0) & = & \frac{3}{16\sqrt{2}} \left(  x_1 ( x_1
x_2 - x_4 x_5 )
( x_1^2 + 7 x_2^2 + 4 x_3^2 +  x_4^2 + 7 x_5^2 + 4 x_6^2))  \right. \\
 &  & - \left. x_3( x_2 x_3 -  x_5 x_6 ) (4 x_1^2 + 7 x_2^2 + x_3^2 + 4 x_4^2 +
7 x_5^2 +  x_6^2 ) \right),\\

\\

 \mathcal{J}_{55}^3(0) & = &  \frac{9}{8\sqrt{2}} \ x_2 ( x_1^2 -
x_3^2 + x_4^2 - x_6^2)( x_1 x_3 - x_4 x_6),   \\

\mathcal{J}_{56}^3(0) & = & \frac{3}{16\sqrt{2}} \left( x_1 (-
x_2^2 (x_5^2 + 4 x_6^2) + x_3^2 (4 x_5^2 + x_6^2) - (x_2^2 -
x_3^2) ( 7 x_1^2 + x_2^2 +
x_3^2 + 7 x_4^2))  \right. \\
 &  & + \left. x_5 x_2 x_4 (7 x_1^2 + x_2^2 + 4 x_3^2 + 7 x_4^2 + x_5^2 + 4 x_6^2) -
 x_6 x_3 x_4 (7 x_1^2 + 4 x_2^2 + x_3^2 + 7 x_4^2 + 4 x_5^2 + x_6^2 )\right), \\

\\

 \mathcal{J}_{66}^3(0) & = &  \frac{9}{8\sqrt{2}} \ x_3 ( x_1^2 -
x_2^2 + x_4^2 - x_5^2)(- x_1 x_2 + x_4 x_5),   \\

\end{array}
\end{equation}

\newpage

$\mathcal{J}_0^{4)} = \left(\mathcal{J}_{ij}^4(0)\right)$, $i, j =
1, \ldots, 6$, where
\[
\begin{array}{ccl}

\mathcal{J}_{11}^4(0) & = & \frac{3}{32} \left( - 8 (x_1^2 + x_2^2
+  x_3^2 + x_4^2 +  x_5^2 +  x_6^2) ( x_2  x_6
 + x_3 x_5 )^2 \right. \\
 &  & + \left. x_4^2 ( (x_1^2 + 7 x_2^2 + 7 x_3^2 + x_4^2 + 7 x_5^2 + 7 x_6^2) ( x_2^2 + x_3^2 + x_5^2 +
x_6^2) - 6 (x_2^2 + x_5^2) (
x_3^2 + x_6^2))  \right) \\

\mathcal{J}_{12}^4(0) & = & \frac{3}{64} \left(  8 x_6 (x_1^2 +
x_2^2 +  x_3^2 + x_4^2 +  x_5^2 +  x_6^2) (x_1 x_2  x_6 + x_1 x_3
x_5 + x_2 x_3 x_4)
  \right. \\
 &  & -  x_4 x_5 ( (x_1^2 +  x_2^2 +  x_3^2 +
x_4^2 +  x_5^2 +  x_6^2)^2 - (x_1^2 +  x_2^2 +  x_3^2 + x_4^2 +
x_5^2 + x_6^2) (3
x_3^2 - 5 x_6^2 ) \\
 &  & + 12 \left. ( (
x_1^2 + x_4^2)(x_2^2 + x_5^2) - (x_3^2 + x_6^2)^2))  \right) \\

\mathcal{J}_{13}^4(0) & = & \frac{3}{64} \left(  8 x_5 (x_1^2 +
x_2^2 +  x_3^2 + x_4^2 +  x_5^2 +  x_6^2) (x_1 x_2  x_6 + x_1 x_3
x_5 + x_2 x_3 x_4)
  \right. \\
 &  & -  x_4 x_6 ( (x_1^2 +  x_2^2 +  x_3^2 +
x_4^2 +  x_5^2 +  x_6^2)^2 + (x_1^2 +  x_2^2 +  x_3^2 + x_4^2 +
x_5^2 + x_6^2) (-3
x_2^2 + 5 x_5^2 ) \\
 &  & + 12 \left. ( (
x_1^2 + x_4^2)(x_3^2 + x_6^2) - (x_2^2 + x_5^2)^2))  \right) \\

\mathcal{J}_{14}^4(0) & = & \frac{3}{32} \left(  8  (x_1^2 + x_2^2
+  x_3^2 + x_4^2 +  x_5^2 +  x_6^2) (- x_2 x_3  + x_5 x_6)( x_2
x_6 + x_3 x_5  )
  \right. \\
 &  & - \left. x_1 x_4 ( (x_1^2 + x_4^2 ) (x_2^2 +  x_3^2 +
  x_5^2 +  x_6^2) - 6  (
x_2^2 + x_5^2)(x_3^2 + x_6^2) + 7 (x_2^2 +  x_3^2 +
  x_5^2 +  x_6^2)^2 \right) \\

  \mathcal{J}_{15}^4(0) & = & \frac{3}{64} \left(  8 x_3 (x_1^2 +  x_2^2 +  x_3^2 +
x_4^2 +  x_5^2 +  x_6^2) (x_1 x_2  x_6 + x_1 x_3 x_5 - x_4 x_5
x_6)
  \right. \\
 &  & +  x_2 x_4 ( (x_1^2 +  x_2^2 +  x_3^2 +
x_4^2 +  x_5^2 +  x_6^2)^2 + (x_1^2 +  x_2^2 +  x_3^2 + x_4^2 +
x_5^2 + x_6^2) (5
x_3^2 - 3 x_6^2 ) \\
 &  & + 12 \left. ( (
x_1^2 + x_4^2)(x_2^2 + x_5^2) - (x_3^2 + x_6^2)^2))  \right) \\

  \mathcal{J}_{16}^4(0) & = & \frac{3}{64} \left(  8 x_2 (x_1^2 +  x_2^2 +  x_3^2 +
x_4^2 +  x_5^2 +  x_6^2) (x_1 x_2  x_6 + x_1 x_3 x_5 - x_4 x_5
x_6)
  \right. \\
 &  & +  x_3 x_4 ( (x_1^2 +  x_2^2 +  x_3^2 +
x_4^2 +  x_5^2 +  x_6^2)^2 + (x_1^2 +  x_2^2 +  x_3^2 + x_4^2 +
x_5^2 + x_6^2) (5
x_2^2 - 3 x_5^2 ) \\
 &  & + 12 \left. ( (
x_1^2 + x_4^2)(x_3^2 + x_6^2) - (x_2^2 + x_5^2)^2))  \right) \\

\\
\mathcal{J}_{22}^4(0) & = & \frac{3}{32} \left( - 8 (x_1^2 + x_2^2
+  x_3^2 + x_4^2 +  x_5^2 +  x_6^2) ( x_1  x_6
 + x_3 x_4 )^2 \right. \\
 &  & + \left. x_5^2 ( (7 x_1^2 +  x_2^2 + 7 x_3^2 + 7 x_4^2 +  x_5^2 + 7 x_6^2) ( x_1^2 + x_3^2 + x_4^2 +
x_6^2) - 6 (x_1^2 + x_4^2) (
x_3^2 + x_6^2))  \right) \\

\mathcal{J}_{23}^4(0) & = & \frac{3}{64} \left(  8 x_4 (x_1^2 +
x_2^2 +  x_3^2 + x_4^2 +  x_5^2 +  x_6^2) (x_1 x_2  x_6 + x_1 x_3
x_5 + x_2 x_3 x_4)
  \right. \\
 &  & -  x_5 x_6 ( (x_1^2 +  x_2^2 +  x_3^2 +
x_4^2 +  x_5^2 +  x_6^2)^2 - (x_1^2 +  x_2^2 +  x_3^2 + x_4^2 +
x_5^2 + x_6^2) (3
x_1^2 - 5 x_4^2 ) \\
 &  & + 12 \left. ( (
x_2^2 + x_5^2)(x_3^2 + x_6^2) - (x_1^2 + x_4^2)^2))  \right) \\

\mathcal{J}_{24}^4(0) & = & \frac{3}{64} \left(  8 x_3 (x_1^2 +
x_2^2 +  x_3^2 + x_4^2 +  x_5^2 +  x_6^2) (x_1 x_2  x_6 + x_2 x_3
x_4 - x_4 x_5 x_6)
  \right. \\
 &  & +  x_1 x_5 ( (x_1^2 +  x_2^2 +  x_3^2 +
x_4^2 +  x_5^2 +  x_6^2)^2 + (x_1^2 +  x_2^2 +  x_3^2 + x_4^2 +
x_5^2 + x_6^2) (5
x_3^2 - 3 x_6^2 ) \\
 &  & + 12 \left. ( (
x_1^2 + x_4^2)(x_2^2 + x_5^2) - (x_3^2 + x_6^2)^2))  \right) \\

\mathcal{J}_{25}^4(0) & = & \frac{3}{32} \left(  8  (x_1^2 + x_2^2
+  x_3^2 + x_4^2 + x_5^2 +  x_6^2) (- x_1 x_3  + x_4 x_6)( x_1 x_6
+ x_3 x_4)
  \right. \\
 &  & - \left. x_2 x_5 ( (x_2^2 + x_5^2 ) (x_1^2 +  x_3^2 +
  x_4^2 +  x_6^2) - 6  (
x_1^2 + x_4^2)(x_3^2 + x_6^2) + 7 (x_1^2 +  x_3^2 +
  x_4^2 +  x_6^2)^2 \right) \\

\mathcal{J}_{26}^4(0) & = & \frac{3}{64} \left(  8 x_1 (x_1^2 +
x_2^2 +  x_3^2 + x_4^2 +  x_5^2 +  x_6^2) (x_1 x_2  x_6 + x_2 x_3
x_4 - x_4 x_5 x_6)
  \right. \\
 &  & +  x_3 x_5 ( (x_1^2 +  x_2^2 +  x_3^2 +
x_4^2 +  x_5^2 +  x_6^2)^2 + (x_1^2 +  x_2^2 +  x_3^2 + x_4^2 +
x_5^2 + x_6^2) (5
x_1^2 - 3 x_4^2 ) \\
 &  & + 12 \left. ( (
x_2^2 + x_5^2)(x_3^2 + x_6^2) - (x_1^2 + x_4^2)^2))  \right) \\

\\

\mathcal{J}_{33}^4(0) & = & \frac{3}{32} \left( 8 (x_1^2 +  x_2^2
+  x_3^2 + x_4^2 + x_5^2 +  x_6^2) ( x_1  x_5
 + x_2 x_4 )^2 \right. \\
 &  & - \left. x_6^2 ( (7 x_1^2 + 7 x_2^2 +  x_3^2 + 7 x_4^2 + 7 x_5^2 +  x_6^2) ( x_1^2 + x_2^2 + x_4^2 +
x_5^2) - 6 (x_1^2 + x_4^2) (
x_2^2 + x_5^2))  \right) \\

\mathcal{J}_{34}^4(0) & = & \frac{3}{64} \left(  8 x_2 (x_1^2 +
x_2^2 +  x_3^2 + x_4^2 +  x_5^2 +  x_6^2) (x_1 x_3  x_5 + x_2 x_3
x_4 - x_4 x_5 x_6)
  \right. \\
 &  & +  x_1 x_6 ( (x_1^2 +  x_2^2 +  x_3^2 +
x_4^2 +  x_5^2 +  x_6^2)^2 + (x_1^2 +  x_2^2 +  x_3^2 + x_4^2 +
x_5^2 + x_6^2) (5
x_2^2 - 3 x_5^2 ) \\
 &  & + 12 \left. ( (
x_1^2 + x_4^2)(x_3^2 + x_6^2) - (x_2^2 + x_5^2)^2))  \right) \\

\mathcal{J}_{35}^4(0) & = & \frac{3}{64} \left(  8 x_1 (x_1^2 +
x_2^2 +  x_3^2 + x_4^2 +  x_5^2 +  x_6^2) (x_1 x_3  x_5 + x_2 x_3
x_4 - x_4 x_5 x_6)
  \right. \\
 &  & +  x_2 x_6 ( (x_1^2 +  x_2^2 +  x_3^2 +
x_4^2 +  x_5^2 +  x_6^2)^2 + (x_1^2 +  x_2^2 +  x_3^2 + x_4^2 +
x_5^2 + x_6^2) (5
x_1^2 - 3 x_4^2 ) \\
 &  & + 12 \left. ( (
x_2^2 + x_5^2)(x_3^2 + x_6^2) - (x_1^2 + x_4^2)^2))  \right) \\

\mathcal{J}_{36}^4(0) & = & \frac{3}{32} \left(  8  (x_1^2 + x_2^2
+  x_3^2 + x_4^2 + x_5^2 +  x_6^2) (- x_1 x_2  + x_4 x_5)( x_1 x_5
+ x_2 x_4)
  \right. \\
 &  & - \left. x_3 x_6 ( (x_3^2 + x_6^2 ) (x_1^2 +  x_2^2 +
  x_4^2 +  x_5^2) - 6  (
x_1^2 + x_4^2)(x_2^2 + x_5^2) + 7 (x_1^2 +  x_2^2 +
  x_4^2 +  x_5^2)^2 \right) \\

  \end{array}
\]

\begin{equation}\label{J4}
\begin{array}{ccl}

\mathcal{J}_{44}^4(0) & = & \frac{3}{32} \left(- 8 (x_1^2 +  x_2^2
+  x_3^2 + x_4^2 + x_5^2 +  x_6^2) ( x_2 x_3
 - x_5 x_6 )^2 \right. \\
 &  & + \left. x_1^2 ( ( x_1^2 + 7 x_2^2 + 7 x_3^2 +  x_4^2 + 7 x_5^2 + 7 x_6^2) ( x_2^2 + x_3^2 + x_5^2 +
x_6^2) - 6 (x_2^2 + x_5^2) (
x_3^2 + x_6^2))  \right) \\

\mathcal{J}_{45}^4(0) & = & \frac{3}{64} \left(  8 x_6 (x_1^2 +
x_2^2 +  x_3^2 + x_4^2 +  x_5^2 +  x_6^2) (- x_1 x_3  x_5 - x_2
x_3 x_4 + x_4 x_5 x_6)
  \right. \\
 &  & -  x_1 x_2 ( (x_1^2 +  x_2^2 +  x_3^2 +
x_4^2 +  x_5^2 +  x_6^2)^2 + (x_1^2 +  x_2^2 +  x_3^2 + x_4^2 +
x_5^2 + x_6^2) (3
x_3^2 - 5 x_6^2 ) \\
 &  & + 12 \left. ( (
x_1^2 + x_4^2)(x_2^2 + x_5^2) - (x_3^2 + x_6^2)^2))  \right) \\

\mathcal{J}_{46}^4(0) & = & \frac{3}{64} \left(  8 x_5 (x_1^2 +
x_2^2 +  x_3^2 + x_4^2 +  x_5^2 +  x_6^2) (- x_1 x_2  x_6 -  x_2
x_3 x_4 + x_4 x_5 x_6)
  \right. \\
 &  & -  x_1 x_3 ( (x_1^2 +  x_2^2 +  x_3^2 +
x_4^2 +  x_5^2 +  x_6^2)^2 - (x_1^2 +  x_2^2 +  x_3^2 + x_4^2 +
x_5^2 + x_6^2) (3
x_2^2 - 5 x_5^2 ) \\
 &  & + 12 \left. ( (
x_1^2 + x_4^2)(x_3^2 + x_6^2) - (x_2^2 + x_5^2)^2))  \right) \\

  \\

\mathcal{J}_{55}^4(0) & = & \frac{3}{32} \left(- 8 (x_1^2 +  x_2^2
+  x_3^2 + x_4^2 + x_5^2 +  x_6^2) ( x_1 x_3
 - x_4 x_6 )^2 \right. \\
 &  & + \left. x_2^2 ( ( 7 x_1^2 +  x_2^2 + 7 x_3^2 + 7 x_4^2 +  x_5^2 + 7 x_6^2) ( x_1^2 + x_3^2 + x_4^2 +
x_6^2) - 6 (x_1^2 + x_4^2) (
x_3^2 + x_6^2))  \right) \\

\mathcal{J}_{56}^4(0) & = & \frac{3}{64} \left(  8 x_4 (x_1^2 +
x_2^2 +  x_3^2 + x_4^2 +  x_5^2 +  x_6^2) (- x_1 x_2  x_6 - x_1
x_3  x_5  + x_4 x_5 x_6)
  \right. \\
 &  & -  x_2 x_3 ( (x_1^2 +  x_2^2 +  x_3^2 +
x_4^2 +  x_5^2 +  x_6^2)^2 - (x_1^2 +  x_2^2 +  x_3^2 + x_4^2 +
x_5^2 + x_6^2) (3
x_1^2 - 5 x_4^2 ) \\
 &  & + 12 \left. ( (
x_2^2 + x_5^2)(x_3^2 + x_6^2) - (x_1^2 + x_4^2)^2))  \right) \\

 \\

\mathcal{J}_{66}^4(0) & = & \frac{3}{32} \left(- 8 (x_1^2 +  x_2^2
+  x_3^2 + x_4^2 + x_5^2 +  x_6^2) ( x_1 x_2
 - x_4 x_5 )^2 \right. \\
 &  & + \left. x_3^2 ( ( 7 x_1^2 + 7 x_2^2 +  x_3^2 + 7 x_4^2 + 7 x_5^2 + x_6^2) ( x_1^2 + x_2^2 + x_4^2 +
x_5^2) - 6 (x_1^2 + x_4^2) (
x_2^2 + x_5^2))  \right) \\

\end{array}
\end{equation}

\newpage

$\mathcal{J}_0^{5)} = \left(\mathcal{J}_{ij}^5(0)\right)$, $i, j =
1, \ldots, 6$, where
\[
\begin{array}{ccl}

\mathcal{J}_{11}^5(0) & = &   - \frac{45}{64\sqrt{2}} \ x_4
(x_1^2+x_2^2+x_3^2+x_4^2+x_5^2+x_6^2) (x_2^2 -
x_3^2 + x_5^2 - x_6^2)(x_2 x_6 + x_3 x_5),   \\

\mathcal{J}_{12}^5(0) & = & \frac{3}{32 \sqrt{2}}
(x_1^2+x_2^2+x_3^2+x_4^2+x_5^2+x_6^2)^2 \left(x_3 (x_4^2 - x_5^2)
+ x_6 (x_1 x_4 - x_2 x_5) \right)
\\ & & - \frac{15}{128\sqrt{2}}
(x_1^2+x_2^2+x_3^2+x_4^2+x_5^2+x_6^2)\left( x_3 ( x_1^2 (x_4^2 - 4
x_5^2) + x_2^2 (4 x_4^2 - x_5^2) \right. \\
 &  & + \left. (x_4^2 - x_5^2) (7 x_3^2 +
x_4^2 + x_5^2 + 7 x_6^2))  + x_1 x_4 x_6 (x_1^2 + 4 x_2^2 + 7
x_3^2 + x_4^2 + 4 x_5^2 + 7 x_6^2) \right. \\
 &  & - \left. x_2 x_5 x_6
 (4 x_1^2 + x_2^2 + 7 x_3^2 + 4 x_4^2 + x_5^2 + 7 x_6^2 )\right), \\

\mathcal{J}_{13}^5(0) & = & - \frac{3}{32 \sqrt{2}}
(x_1^2+x_2^2+x_3^2+x_4^2+x_5^2+x_6^2)^2 \left(x_2 (x_4^2 - x_6^2)
+ x_5 (x_1 x_4 - x_3 x_6) \right) \\ & & - \frac{15}{128\sqrt{2}}
(x_1^2+x_2^2+x_3^2+x_4^2+x_5^2+x_6^2) \left( - x_4 ( x_1 x_5 + x_2
x_4 ) (x_1^2 + 7 x_2^2 + 4 x_3^2   \right. \\
 &  & + \left. x_4^2 +7 x_5^2 + 4 x_6^2)+
x_6( x_2 x_6 +  x_3 x_5 ) (4 x_1^2 + 7 x_2^2 + x_3^2 + 4 x_4^2 +
7 x_5^2 + x_6^2 ) \right),\\

\mathcal{J}_{14}^5(0) & = & - \frac{45}{128\sqrt{2}} \
(x_1^2+x_2^2+x_3^2+x_4^2+x_5^2+x_6^2)(x_2^2 - x_3^2 + x_5^2 -
x_6^2)(x_1 (x_2 x_6 - x_3 x_5)   \\
 &  & +  x_4 (x_2 x_3 - x_5 x_6)),  \\

\mathcal{J}_{15}^5(0) & = & \frac{3}{32 \sqrt{2}}
(x_1^2+x_2^2+x_3^2+x_4^2+x_5^2+x_6^2)^2 \left(x_6 (x_2^2 - x_4^2)
+ x_3 (x_1 x_4 + x_2 x_5) \right)
\\ & &- \frac{15}{128\sqrt{2}}
(x_1^2+x_2^2+x_3^2+x_4^2+x_5^2+x_6^2)\left(  x_4 ( x_1 x_3 - x_4
x_6 )(x_1^2 + 4 x_2^2 + 7 x_3^2  \right. \\
 &  & + \left.  x_4^2 +4 x_5^2 + 7 x_6^2))+ x_2( x_2 x_6 +  x_3 x_5 ) (4 x_1^2 + x_2^2 + 7 x_3^2 + 4 x_4^2 + x_5^2 + 7 x_6^2 ) \right),\\

 \mathcal{J}_{16}^5(0) & = & - \frac{3}{32 \sqrt{2}} (x_1^2+x_2^2+x_3^2+x_4^2+x_5^2+x_6^2)^2 \left(x_5
(x_3^2 - x_4^2) +  x_2 (x_1 x_4 + x_3 x_6) \right)
 \\ & & - \frac{15}{128\sqrt{2}}(x_1^2+x_2^2+x_3^2+x_4^2+x_5^2+x_6^2) \left( x_5 ( x_1^2 (x_4^2 - 4
x_3^2) + x_6^2 (4 x_4^2 - x_3^2) \right. \\
 &  & + \left. (x_4^2 - x_3^2) (7 x_2^2 + x_3^2 + x_4^2 + 7 x_5^2)) -  x_1 x_2 x_4 (x_1^2 + 7 x_2^2 + 4 x_3^2 + x_4^2 + 7 x_5^2 + 4 x_6^2)  \right. \\
 &  & - \left. x_6 x_2
 x_3 (4 x_1^2 + 7 x_2^2 + x_3^2 + 4 x_4^2 + 7 x_5^2 + x_6^2 )\right), \\

\\

 \mathcal{J}_{22}^5(0) & = & - \frac{45}{64\sqrt{2}} \ x_5 (x_1^2+x_2^2+x_3^2+x_4^2+x_5^2+x_6^2)( x_1^2 -
x_3^2 + x_4^2 - x_6^2)(x_1 x_6 + x_3 x_4),   \\

\mathcal{J}_{23}^5(0) & = & \frac{3}{32 \sqrt{2}}
(x_1^2+x_2^2+x_3^2+x_4^2+x_5^2+x_6^2)^2\left(x_1 (x_5^2 - x_6^2) +
x_4 (x_2 x_5 - x_3 x_6) \right)\\ & & -\frac{15}{128\sqrt{2}}
(x_1^2+x_2^2+x_3^2+x_4^2+x_5^2+x_6^2)\left(  x_5 ( x_1 x_5 + x_2
x_4 )(7 x_1^2 + x_2^2 + 4 x_3^2   \right. \\
 &  & + \left. 7 x_4^2+x_5^2 + 4 x_6^2))  - x_6( x_1 x_6 +  x_3 x_4 ) (7 x_1^2 +
4 x_2^2 + x_3^2 + 7 x_4^2 +
4 x_5^2 + x_6^2 ) \right),\\

\mathcal{J}_{24}^5(0) & = & - \frac{3}{32 \sqrt{2}}
(x_1^2+x_2^2+x_3^2+x_4^2+x_5^2+x_6^2)^2 \left(x_6 (x_1^2 - x_5^2)
+ x_3 (x_1 x_4 + x_2 x_5) \right) \\ & & -\frac{15}{128\sqrt{2}}
(x_1^2+x_2^2+x_3^2+x_4^2+x_5^2+x_6^2)\left( x_6 ( x_2^2 ( 4 x_1^2
- x_5^2 ) + x_4^2 ( x_1^2 - 4 x_5^2) \right. \\ &  & + \left.
(x_1^2 - x_5^2) (x_1^2 + 7 x_3^2 + x_5^2 + 7 x_6^2))
-x_2 x_3 x_5 (4 x_1^2 + x_2^2 + 7 x_3^2 + 4 x_4^2 + x_5^2 + 7 x_6^2 )  \right. \\
 &  & - \left. x_4 x_1 x_3
  (x_1^2 + 4 x_2^2 + 7 x_3^2 + x_4^2 + 4 x_5^2 + 7 x_6^2)  \right), \\

\mathcal{J}_{25}^5(0) & = &  - \frac{45}{128\sqrt{2}} \
(x_1^2+x_2^2+x_3^2+x_4^2+x_5^2+x_6^2)(x_1^2 - x_3^2 + x_4^2 -
x_6^2)(x_2 (x_1 x_6 + x_3 x_4) \\
 &  & -  x_5 (x_1 x_3 + x_4 x_6)),  \\

\mathcal{J}_{26}^5(0) & = & \frac{3}{32 \sqrt{2}}
(x_1^2+x_2^2+x_3^2+x_4^2+x_5^2+x_6^2)^2 \left(x_4 (x_3^2 - x_5^2)
+  x_1 (x_2 x_5 + x_3 x_6) \right) \\ & & -\frac{15}{128\sqrt{2}}
(x_1^2+x_2^2+x_3^2+x_4^2+x_5^2+x_6^2)\left( x_4 ( x_2^2 (4 x_3^2 -
x_5^2) + x_6^2 (x_3^2 - 4 x_5^2)\right. \\
&  & + \left. (x_3^2 - x_5^2) (7 x_1^2 +x_3^2 + 7 x_4^2 + x_5^2))
+ x_2 x_1 x_5 (7 x_1^2 + x_2^2 + 4 x_3^2 + 7 x_4^2 + x_5^2 + 4
x_6^2) \right. \\
&  & + \left. x_6
 x_1 x_3 (7 x_1^2 + 4 x_2^2 + x_3^2 + 7 x_4^2 + 4 x_5^2 + x_6^2 )\right), \\

\\

 \mathcal{J}_{33}^5(0) & = &  -\frac{45}{64\sqrt{2}} \ x_6 (x_1^2+x_2^2+x_3^2+x_4^2+x_5^2+x_6^2)( x_1^2 -
x_2^2 + x_4^2 - x_5^2)(x_1 x_5 + x_2 x_4),   \\

\mathcal{J}_{34}^5(0) & = & \frac{3}{32 \sqrt{2}}
(x_1^2+x_2^2+x_3^2+x_4^2+x_5^2+x_6^2)^2 \left(x_5 (x_1^2 - x_6^2)
+ x_2 (x_1 x_4 + x_3 x_6) \right) \\ & & -\frac{15}{128\sqrt{2}}
(x_1^2+x_2^2+x_3^2+x_4^2+x_5^2+x_6^2)\left(  x_6 ( x_2 x_3 - x_5
x_6 )(4 x_1^2 + 7 x_2^2 + x_3^2 \right. \\
 &  & + \left. 4 x_4^2 + 7 x_5^2 + x_6^2)) + x_1( x_1 x_5 +  x_2 x_4 ) (x_1^2 + 7 x_2^2 + 4 x_3^2 + x_4^2 +
7 x_5^2 + 4 x_6^2 ) \right),\\

\mathcal{J}_{35}^5(0) & = & - \frac{3}{32 \sqrt{2}}
(x_1^2+x_2^2+x_3^2+x_4^2+x_5^2+x_6^2)^2 \left(x_4 (x_2^2 - x_6^2)
+ x_1 (x_2 x_5 + x_3 x_6) \right)\\ & & -\frac{15}{128\sqrt{2}}
(x_1^2+x_2^2+x_3^2+x_4^2+x_5^2+x_6^2)\left(  x_6 (- x_1 x_3 + x_4
x_6 )(7 x_1^2 + 4 x_2^2 + x_3^2 + 7 x_4^2 \right. \\
 &  & + \left. 4 x_5^2 + x_6^2))
-x_2( x_1 x_5 +  x_2 x_4 ) (7 x_1^2 + x_2^2 + 4 x_3^2 + 7 x_4^2 +
x_5^2 + 4 x_6^2 ) \right),\\

\mathcal{J}_{36}^5(0) & = &  - \frac{45}{128\sqrt{2}} \
(x_1^2+x_2^2+x_3^2+x_4^2+x_5^2+x_6^2)(x_1^2 - x_2^2 + x_4^2 -
x_5^2)(- x_3 (x_1 x_5 + x_2 x_4) \\
 &  & + x_6 (x_1 x_2 - x_4 x_5)),  \\

\end{array}
\]

\begin{equation}\label{J5}
\begin{array}{ccl}

%\\
 \mathcal{J}_{44}^5(0) & = &  - \frac{45}{64\sqrt{2}} \ x_1 (x_1^2+x_2^2+x_3^2+x_4^2+x_5^2+x_6^2)( x_2^2 -
x_3^2 + x_5^2 - x_6^2)(- x_2 x_3 + x_5 x_6),   \\

\mathcal{J}_{45}^5(0) & = & \frac{3}{32 \sqrt{2}}
(x_1^2+x_2^2+x_3^2+x_4^2+x_5^2+x_6^2)^2\left(x_3 (x_2^2 - x_1^2) +
x_6 (x_1 x_4 - x_2 x_5) \right)
\\ &  & -\frac{15}{128\sqrt{2}}
(x_1^2+x_2^2+x_3^2+x_4^2+x_5^2+x_6^2)\left( x_3 (- x_1^2 (x_4^2 +
4 x_5^2) + x_2^2 (4 x_4^2 + x_5^2)  \right. \\
 &  & - \left. (x_1^2 - x_2^2) ( x_1^2 +
x_2^2 + 7 x_3^2 + 7 x_6^2)) + x_4 x_1 x_6 (x_1^2 + 4 x_2^2 + 7
x_3^2 + x_4^2 + 4 x_5^2 + 7 x_6^2) \right. \\
 &  & + \left. x_5 x_2 x_6 (4 x_1^2 + x_2^2 + 7 x_3^2 + 4 x_4^2 +  x_5^2 + 7 x_6^2 )\right), \\

\mathcal{J}_{46}^5(0) & = & -\frac{3}{32 \sqrt{2}}
(x_1^2+x_2^2+x_3^2+x_4^2+x_5^2+x_6^2)^2\left(x_2 (x_3^2 - x_1^2) +
x_5 (x_1 x_4 - x_3 x_6) \right)\\ &  & -\frac{15}{128\sqrt{2}}
(x_1^2+x_2^2+x_3^2+x_4^2+x_5^2+x_6^2)\left(  x_1 ( x_1 x_2 - x_4
x_5 )( x_1^2 + 7 x_2^2 + 4 x_3^2 +  x_4^2 \right. \\
 &  & + \left. 7 x_5^2 + 4 x_6^2))
- x_3( x_2 x_3 -  x_5 x_6 ) (4 x_1^2 + 7 x_2^2 + x_3^2 + 4 x_4^2 +
7 x_5^2 +  x_6^2 ) \right),\\

\\

 \mathcal{J}_{55}^5(0) & = & - \frac{45}{64\sqrt{2}} \ x_2 (x_1^2+x_2^2+x_3^2+x_4^2+x_5^2+x_6^2)( x_1^2 -
x_3^2 + x_4^2 - x_6^2)( x_1 x_3 - x_4 x_6),   \\

\mathcal{J}_{56}^5(0) & = & \frac{3}{32 \sqrt{2}}
(x_1^2+x_2^2+x_3^2+x_4^2+x_5^2+x_6^2)^2 \left(x_1
(x_3^2 - x_2^2) +  x_4 (x_2 x_5 - x_3 x_6) \right) \\
 &  & -\frac{15}{128\sqrt{2}}
(x_1^2+x_2^2+x_3^2+x_4^2+x_5^2+x_6^2)\left( x_1 (- x_2^2 (x_5^2 +
4 x_6^2) + x_3^2 (4 x_5^2 + x_6^2) \right. \\
 &  & - \left. (x_2^2 - x_3^2) ( 7 x_1^2 +
x_2^2 + x_3^2 + 7 x_4^2)) + x_5 x_2 x_4 (7 x_1^2 + x_2^2 + 4 x_3^2
+ 7 x_4^2 + x_5^2 + 4 x_6^2)  \right. \\
 &  & - \left.
 x_6 x_3 x_4 (7 x_1^2 + 4 x_2^2 + x_3^2 + 7 x_4^2 + 4 x_5^2 + x_6^2 )\right), \\

\\

 \mathcal{J}_{66}^5(0) & = &  - \frac{45}{64\sqrt{2}} \ x_3 (x_1^2+x_2^2+x_3^2+x_4^2+x_5^2+x_6^2)( x_1^2 -
x_2^2 + x_4^2 - x_5^2)(- x_1 x_2 + x_4 x_5),   \\

 \end{array}
\end{equation}

\end{document}